# Forecasting unstable processes

**Jin-Lung Lin[1] and Ching-Zong Wei[2]**

*Academia Sinica*

**Abstract:** Previous analysis on forecasting theory either assume knowing the true parameters or assume the stationarity of the series. Not much are known on the forecasting theory for nonstationary process with estimated parameters. This paper investigates the recursive least square forecast for stationary and nonstationary processes with unit roots. We first prove that the accumulated forecast mean square error can be decomposed into two components, one of which arises from estimation uncertainty and the other from the disturbance term. The former, of the order of $\log(T)$, is of second order importance to the latter term, of the order T. However, since the latter is common for all predictors, it is the former that determines the property of each predictor. Our theorem implies that the improvement of forecasting precision is of the order of $\log(T)$ when existence of unit root is properly detected and taken into account. Also, our theorem leads to a new proof of strong consistency of predictive least squares in model selection and a new test of unit root where no regression is needed.

The simulation results confirm our theoretical findings. In addition, we find that while mis-specification of AR order and under-specification of the number of unit root have marginal impact on forecasting precision, over-specification of the number of unit root strongly deteriorates the quality of long term forecast. As for the empirical study using Taiwanese data, the results are mixed. Adaptive forecast and imposing unit root improve forecast precision for some cases but deteriorate forecasting precision for other cases.

## 1. Introduction

Forecasting future observations is one of the major purpose of building a time series model. Even for the purpose of time series controlling, forecasting provide the essential basis. For this purpose, autoregressive (AR) models are widely used for their simplicity. For an AR(p) process,

$$y_t = \beta_1 y_{t-1} + \beta_2 y_{t-2} + \cdots + \beta_p y_{t-p} + \epsilon_t \tag{1}$$

where $\phi(z) = 1 - \beta_1 z - \cdots - \beta_p z^p$ the characteristic polynomial determines the properties of the series. $y_t$ is called stationary or stable if all roots of $\phi$ are outside the unit circle, unstable or nonstationary if some roots of $\phi$ are on the unit circle and explosive if some roots of $\phi$ are inside the unit circle. Previous analysis on forecasting theory either assume knowing true $\beta_s'$ or only consider the stationary cases. For examples, Ing [8, 9] and Bhansali [1, 2] analyze the multistep prediction of stationary AR processes while Ing [7] derives the mean squares prediction errors of the least squares predictors in random walk model. Not much are known on the forecasting theory for unstable process with estimated parameters. This paper investigates the recursive least square forecast for stable and unstable processes.

[1]Institute of Economics, Academia Sinica, 128 Sec. 2, Academia Rd., Nankang, Taipei, Taiwan 11529, e-mail: jlin@econ.sinica.edu.tw

[2]Institute of Statistical Science, Academia Sinica, 128 Sec. 2, Academia Rd., Nankang, Taipei, Taiwan 11529.

Let $\hat{y}_t$ be the forecast of $y_t$ based upon information up to $t-1$. If one is interested in one-period forecast, $(y_t - \hat{y}_t)^2$ is the cost to be minimized. However, there are two situations where the accumulated cost function, $\sum_{k=1}^{t}(y_k - \hat{y}_k)^2$ is more appropriate. First, in the sequential forecast case, (see Goodwin and Sin [6]) the forecaster are updated sequentially over many periods and the accumulated cost function is the target to be minimized. Second, for a single realization of time series, the averaged accumulated cost function is often used as the yardstick to evaluate the out-of-sample forecasting performance of alternative forecasters.

Ing [7] advocated adopting the accumulated cost function $\sum_{t=1}^{T} E(y_t - \hat{y}_t)^2$ over the one-period expected loss function $E(y_{T+1} - \hat{y}_{T+1})^2$. For an AR(1) process, these two quantities are respectively:

$$\frac{1}{T-2}\sum_{t=3}^{T} E(y_t - \hat{y}_t)^2 = \sigma^2 + \frac{2\sigma^2 \log(T)}{T} + o\left(\frac{\log(T)}{T}\right)$$

$$E(y_{T+1} - \hat{y}_{T+1})^2 = \sigma^2 + \frac{2\sigma^2}{T} + o\left(\frac{1}{T}\right)$$

when true $\beta_1 = 1$. In other words, the efficiency loss for not taking the unit root into consideration is greater for the accumulated cost function than the one-period cost function. See also Ing and Wei [11]. It is worth mentioning that Rissanen [14] predictive least square (PLS) for model selection built upon accumulated cost function minimization. See also Wei [18].

Under the assumption that $E(\epsilon_t^2|\mathcal{F}_{t-1}) = \sigma^2$ a.s. for all t, where $\mathcal{F}_{t-1}$ is the sigma field generated by $\{x_s, s \leq t-1\}$, then it can be shown that under appropriate assumptions that $\frac{1}{T}\sum_{t=1}^{T}(y_t - \hat{y}_t)^2 \longrightarrow \sigma^2 \quad a.s.$ But by Chow [4], it is seen that

$$\sum_{t=1}^{T}(y_t - \hat{y}_t)^2 = \sum_{t=1}^{T}\epsilon_t^2 + C_T(1 + o(1)) \quad a.s. \quad \text{on the set } \{C_T \to \infty\}$$

$$\sum_{t=1}^{T}(y_t - \hat{y}_t)^2 = \sum_{t=1}^{T}\epsilon_t^2 + C_T(1 + O(1)) \quad a.s. \quad \text{on the set } \{C_T < \infty\}$$

where

$$C_T = \sum_{t=1}^{T}(y_t - \hat{y}_t - \epsilon_t)^2$$

While $\sum_{t=1}^{T}\epsilon_t^2$ is larger in order than $C_T$, it is common for all forecasters and cannot be removed. Hence $C_T$ becomes a more important quantity when evaluating the performance of alternative forecasters.

Let $\hat{\beta}_t$ be the least square estimate of $\beta$

$$\hat{\beta}_t = [\sum_{k=1}^{t}\boldsymbol{Y}_{k-1}\boldsymbol{Y}'_{k-1}]^{-1}\sum_{k=1}^{t}\boldsymbol{Y}'_{k-1}y_k$$

where $Y_t = \{y_1, \ldots, y_t\}'$, then $\hat{y}_t = \hat{\beta}'_{t-1}Y_{t-1}$ is the least square prediction of $y_t$ at time $t-1$.
Let

$$\phi(z) = (z-1)^a(z+1)^b\Pi_{k=1}^{l}(z^2 - 2\cos\theta_k z + 1)^{d_k}\pi(z)$$

where all roots of $\pi(z)$ are all outside the unit circle. Wei [17] proves that,

$$C_T \to (p + a^2 + b^2 + 2\sum_{k=1}^{l} d_k^2)\sigma^2 \log(T) \quad \text{in probability.} \tag{2}$$

In other words, when $\phi(z)$ has multiple unit roots the accumulated loss increase not linearly with the number of unit roots but at the rate of the square of the number of unit roots.

In this paper, we prove that when $\phi(z)$ has no complex roots, the convergence in (2) can be improved to be almost surely. This result could lead to a new proof of strong consistency of PLS in AR model selection. It is also conjectured that the result of almost surely convergence hold for the case of complex unit roots. We conduct several simulation experiments to assess the convergence result for various sample sizes. In addition, we also consider the impact of near unit root and model mis-specification on multi-step forecasting. Finally, we apply our methods to six real macroeconomic series in Taiwan. Forecasting performance of various forecasters and adaptive forecaster are investigated.

The rest of the paper is organized as follows. The proof of the main theorem is put in Section 2. Section 3 illustrates implications and applications of our main theorem. Section 4 discusses multi-step and adaptive forecast. Monte Carlo results are reported in Section 5 and Section 6 summarizes the empirical results. Section 7 concludes.

## 2. Main theorem

Assume that $\epsilon_t$ are *i.i.d.* random variables with $E(\epsilon_t) = 0$ and $0 < E(\epsilon_t^2) = \sigma^2 < \infty$. Let $\boldsymbol{X}_t = (x_{t-1}, \ldots, x_{t-p})', S_T = \sum_{t=1}^T \epsilon_t$ and $T_T = (-1)^T \sum_{t=1}^T (-1)^t \epsilon_t = \epsilon_t + (-1)T_{T-1}$.

**Lemma 1.** *Assume that* $\boldsymbol{X}_{t+1} = A\boldsymbol{X}_t + \boldsymbol{\varepsilon}_t$ *, where* $\boldsymbol{\epsilon}_t = (\epsilon_t, 0, \ldots, 0)'$ *and the eigenvalues of* $A$ *are all inside the unit circle. Then*

$$\lim_{T\to\infty} \frac{\sum_{t=1}^T \boldsymbol{X}_t S_t}{\sqrt{T \sum_{t=1}^T S_t^2}} = 0 \quad a.s.$$

*and*

$$\lim_{T\to\infty} \frac{\sum_{t=1}^T \boldsymbol{X}_t T_t}{\sqrt{T \sum_{t=1}^T T_t^2}} = 0 \quad a.s.$$

*Proof.* It is known from Lai and Wei [12] [pages 363 and 364] that

$$\lim_{T\to\infty} \frac{1}{T} \sum_{t=1}^T \boldsymbol{X}_t \boldsymbol{X}_t' = \Sigma \quad a.s. \tag{3}$$

where $\Sigma$ is a positive definite matrix,

$$\limsup_{T\to\infty} \frac{\sum_{t=1}^T S_t^2}{T^2 \log\log(T)} = \frac{8\sigma^2}{\pi^2} \quad a.s. \tag{4}$$

and

$$\liminf_{T\to\infty} \frac{\sum_{t=1}^{T} S_t^2}{T^2/\log\log(T)} = \frac{\sigma^2}{4} \quad a.s. \tag{5}$$

Let $\|u\|$ denote the Euclidean norm of a $k$-dimensional vector $u = (u_1, \ldots, u_k)'$, i.e., $\|u\|^2 = \sum_{i=1}^{k} u_i^2$. By (3), $\frac{\|\boldsymbol{X}_T\|^2}{T} \to 0 \quad a.s.$ and in turn we have that

$$0 \le \boldsymbol{X}'_T(\sum_{t=1}^{T} \boldsymbol{X}_t \boldsymbol{X}'_t)^{-1} \boldsymbol{X}_T \le \frac{\|\boldsymbol{X}_T\|^2}{\lambda_{\min}(\sum_{t=1}^{T} \boldsymbol{X}_t \boldsymbol{X}'_t)} = \frac{\|\boldsymbol{X}_T\|^2/T}{\lambda_{\min}(\frac{1}{T}\sum_{t=1}^{T} \boldsymbol{X}_t \boldsymbol{X}'_t)} \to 0 \ a.s.$$

and

$$\boldsymbol{X}'_T(\sum_{t=1}^{T} \boldsymbol{X}_t \boldsymbol{X}'_t)^{-1} \boldsymbol{X}_T \to 0 \quad a.s. \tag{6}$$

where $\lambda_{\min}(A)$ denotes the minimal eigenvalue of matrix $A$.

Furthermore, by the law of iterative logarithm,

$$\limsup_{T\to\infty} \frac{S_T^2}{2T\log\log T} = \sigma^2 \quad a.s.$$

Hence (5) implies that

$$\begin{aligned} \frac{S_T^2}{\sum_{t=1}^{T} S_t^2} &= O\left(\frac{T\log\log T}{T^2/\log\log T}\right) \\ &= O\left(\frac{(\log\log T)^2}{T}\right) \\ &= o(1) \quad a.s. \end{aligned} \tag{7}$$

Now, let

$$\boldsymbol{Z}_T = \frac{\sum_{t=1}^{T} \boldsymbol{X}_t S_t}{(T\sum_{t=1}^{T} S_t^2)^{1/2}}.$$

Then

$$\begin{aligned} \boldsymbol{Z}_T - \boldsymbol{Z}_{T-1} &= \boldsymbol{Z}_T - \frac{\sum_{t=1}^{T-1} \boldsymbol{X}_t S_t}{(T\sum_{t=1}^{T} S_t^2)^{1/2}} - \boldsymbol{Z}_{T-1}\left(1 - (\frac{(T-1)\sum_{t=1}^{T-1} S_t^2}{T\sum_{t=1}^{T} S_t^2})^{1/2}\right) \\ &= \frac{\boldsymbol{X}_T S_T}{(T\sum_{t=1}^{T} S_t^2)^{1/2}} - \boldsymbol{Z}_{T-1}\left(1 - (\frac{T-1}{T} - \frac{T-1}{T}\frac{S_T^2}{\sum_{t=1}^{T} S_t^2})^{1/2}\right) \\ &= \frac{\boldsymbol{X}_T S_T}{(T\sum_{t=1}^{T} S_t^2)^{1/2}} - \boldsymbol{Z}_{T-1}\, o(1), \quad \text{by (7)} \\ &= o(1) - o(1), \quad \text{since} \quad \sup_T \|\boldsymbol{Z}_T\| \le \{\frac{1}{T}\sum_{t=1}^{T} \|\boldsymbol{X}_t\|^2\}^{1/2} \quad a.s. \\ &= o(1) \end{aligned} \tag{8}$$

But,

$$\begin{aligned}
\sum_{t=1}^{T} \boldsymbol{X}_t S_t &= \sum_{t=1}^{T} (A\boldsymbol{X}_{t-1} + \boldsymbol{\varepsilon}_t) S_t \\
&= A\sum_{t=1}^{T} \boldsymbol{X}_{t-1} S_{t-1} + A\sum_{t=1}^{T} \boldsymbol{X}_{t-1}\epsilon_t + \sum_{t=1}^{T} \boldsymbol{\varepsilon}_t S_{t-1} + \sum_{t=1}^{T} \epsilon_t^2 \\
&= A\sum_{t=1}^{T} \boldsymbol{X}_{t-1} S_{t-1} + o((\sum_{t=1}^{T} \|\boldsymbol{X}_{t-1}\|^2)^{1/2} (\log \sum_{t=1}^{T} \|\boldsymbol{X}_{t-1}\|^2)^{\frac{1+\sigma}{2}}) \\
&\quad + o((\sum_{t=1}^{T} S_{t-1}^2)^{1/2} (\log \sum_{t=1}^{T} S_{t-1}^2)^{\frac{1+\sigma}{2}}) + O(T) \quad a.s. \\
&= A(\sum_{t=1}^{T} \boldsymbol{X}_{t-1} S_{t-1}) + o(T^{1/2} (\log T)^{\frac{1+\sigma}{2}}) \\
&\quad + o((\sum_{t=1}^{T} S_{t-1}^2)^{1/2} (\log T)^{\frac{1+\sigma}{2}}) + O(T)
\end{aligned}$$

This implies that

$$\boldsymbol{Z}_T = A\boldsymbol{Z}_{T-1}(1 + o(1)) + o(1) \quad a.s. \tag{9}$$

Combining (8) and (9), we have that

$$\boldsymbol{Z}_{T-1} - A\boldsymbol{Z}_{T-1} = o(1) \quad a.s. \tag{10}$$

Therefore, any limit point $z$ of $\{z_T\}$ would satisfy

$$\boldsymbol{Z} - A\boldsymbol{Z} = 0 \tag{11}$$

Since 1 is not an eigenvalue of A, $\boldsymbol{Z} = 0$. Using the same method one can prove that

$$\frac{\sum_{t=1}^{T} \boldsymbol{X}_t T_t}{(T\sum_{t=1}^{T} T_t^2)^{1/2}} = 0 \quad a.s.$$

This proves Lemma 1. □

**Lemma 2.** *If $E|\epsilon_t^\alpha| < \infty$ for some $\alpha > 2$, then*

$$\lim_{T\to\infty} \frac{\sum_{t=1}^{T} S_t T_t}{\sqrt{\sum_{t=1}^{T} S_t^2 \sum_{t=1}^{T} T_t^2}} = 0 \quad a.s.$$

*Proof.* Note that $\tilde{T}_T = (-1)^T T_T = \sum_{t=1}^{T} (-1)^T \epsilon_t$. Using theorem 3.2 of Phillip in page 234 of Eberlein and Taqqu [5], (4) and (5) hold if we replace $S_t$ by $T_t$.

Therefore,

$$\frac{T_T^2}{\sum_{t=1}^{T} T_t^2} = \frac{\tilde{T}_T^2}{\sum_{t=1}^{T} \tilde{T}_t^2} \to 0 \quad a.s.$$

Let

$$u_T = \frac{\sum_{t=1}^T S_t T_t}{\sqrt{\sum_{t=1}^T S_t^2 \sum_{t=1}^T T_t^2}}.$$

Then

$$\begin{aligned} u_T - u_{T-1} &= \frac{S_T T_T}{\sqrt{\sum_{t=1}^T S_t^2 \sum_{t=1}^T T_t^2}} + u_{T-1}\left(\sqrt{\frac{\sum_{t=1}^{T-1} S_t^2 \sum_{t=1}^T T_t^2}{\sum_{t=1}^T S_t^2 \sum_{t=1}^T T_t^2}} - 1\right) \\ &= o(1) \quad a.s. \end{aligned} \tag{12}$$

But

$$\begin{aligned} \sum_{t=1}^T S_t T_t &= \sum_{t=1}^T (S_{t-1} + \epsilon_t)(-T_{t-1} + \epsilon_t) \\ &= -\sum_{t=1}^T S_{t-1}T_{t-1} + \sum_{t=1}^T S_{t-1}\epsilon_t - \sum_{t=1}^T T_{t-1}\epsilon_t + \sum_{t=1}^T \epsilon_t^2 \\ &= -\sum_{t=1}^{T-1} S_t T_t + o((\sum_{t=1}^T S_{t-1}^2)^{1/2})(\log(\sum_{t=1}^T S_{t-1}^2)) \\ &\quad + o((\sum_{t=1}^T T_t^2)^{1/2}(\log(\sum_{t=1}^T T_t^2))) + O(T) \quad a.s. \end{aligned}$$

Therefore,

$$\begin{aligned} u_T &= -\frac{\sum_{t=1}^{T-1} S_t T_t}{\sqrt{\sum_{t=1}^T S_t^2}\sqrt{\sum_{t=1}^T T_t^2}} + o\left(\frac{\log(\sum_{t=1}^T S_{t-1}^2)}{\sqrt{\sum_{t=1}^T T_t^2}}\right) + o\left(\frac{\log(\sum_{t=1}^T T_t^2)}{\sqrt{\sum_{t=1}^T S_t^2}}\right) \\ &\quad + o\left(\frac{T}{\sqrt{\sum_{t=1}^T T_t^2 \sum_{t=1}^T S_t^2}}\right) \\ &= -U_{T-1}(1 + o(1)) + o(1) \quad a.s. \\ &= -u_{T-1} + o(1) \quad a.s. \end{aligned} \tag{13}$$

Combining (12) and (13), since

$$u_T = o(1) \quad a.s. \quad u_T \longrightarrow 0 \quad a.s.$$ □

Now, we are ready to state our main result.

Let

$$y_t = \beta_1 y_{t-1} + \cdots + \beta_p y_{t-p} + \epsilon_t \tag{14}$$

be an AR(p) model with

$$\phi(z) = 1 - \beta_1 z - \cdots - \beta_p z^p \tag{15}$$

$$= (1 - z)(1 + z)\Psi(z) \tag{16}$$

where $\Psi(z) = 1 - \Psi_1 z - \cdots - \Psi_q z^q$ is a polynomial of order $q = p - 2$ which has all roots outside the unit circle.

**Theorem 1.** *Assume that the AR(p) model (14) satisfies (16). If $\{\epsilon_t\}$ is a sequence of i.i.d. random variables with $E|\epsilon_t|^\alpha < \infty$, where $\alpha > 2$, and $y_0, \ldots, y_{1-p}$ is independent of $\{\epsilon_t\}$ then*

$$\lim_{T\to\infty} \frac{1}{\log T} \log\det(\sum_{t=1}^{T} \boldsymbol{y}_t\boldsymbol{y}'_t) = (p+2) \quad a.s. \tag{17}$$

*where $\boldsymbol{y}'_t = (y_t, \ldots, y_{t-p+1})$.*

*Proof.* By Chan and Wei [3] there exists a non-singular $p \times p$ matrix $Q$ such that $Q\boldsymbol{y}_t = (u_t, v_t, \boldsymbol{x}'_t)$, where

$$\begin{aligned}
\boldsymbol{x}_t &= (x_{t-1}, \ldots, x_{t-q})', \\
u_t &= u_{t-1} + \epsilon_t, \\
v_t &= -v_{t-1} + \epsilon_t \quad \text{and} \\
x_t &= \Psi_1 x_{t-1} + \cdots + \Psi_q x_{t-q}.
\end{aligned}$$

Therefore, if we let $\boldsymbol{z}_t = Q\boldsymbol{y}_t$,

$$\det(\sum_{t=1}^{T} \boldsymbol{y}_t\boldsymbol{y}'_t) = \det[Q^{-1}\sum_{t=1}^{T} \boldsymbol{z}_t\boldsymbol{z}'_t Q^{-1}] = \frac{\det(\sum_{t=1}^{T} \boldsymbol{z}_t\boldsymbol{z}'_t)}{(\det(Q))^2}.$$

To show (17), it is sufficient to show

$$\lim_{T\to\infty} \frac{1}{\log T} \log\det(\sum_{t=1}^{T} \boldsymbol{z}_t\boldsymbol{z}'_t) = (p+2) \quad a.s. \tag{18}$$

Let

$$G_T = \begin{pmatrix} (\sum_{t=1}^{T} u_t^2)^{-1/2} & 0 & 0 \\ 0 & (\sum_{t=1}^{T} v_t^2)^{-1/2} & 0 \\ 0 & 0 & T^{-1/2} I_q \end{pmatrix},$$

where $I_q$ is the $q \times q$ identity matrix.

Then

$$G_T \sum_{t=1}^{T} \boldsymbol{z}_t\boldsymbol{z}'_t G_T = \begin{pmatrix} 1 & a_T & \boldsymbol{b}'_T \\ a_T & 1 & \boldsymbol{c}'_T \\ \boldsymbol{b}_T & \boldsymbol{c}_T & \Gamma_T \end{pmatrix},$$

where

$$\begin{aligned}
a_T &= \frac{(\sum_{t=1}^{T} u_t v_t)}{[(\sum_{t=1}^{T} u_t^2)(\sum_{t=1}^{T} v_t^2)]^{1/2}}, \\
\boldsymbol{b}_T &= \frac{\sum_{t=1}^{T} u_t \boldsymbol{x}_t}{(T\sum_{t=1}^{T} u_t^2)^{1/2}}, \\
\boldsymbol{c}_T &= \frac{\sum_{t=1}^{T} v_t \boldsymbol{x}_t}{(T\sum_{t=1}^{T} v_t^2)^{1/2}},
\end{aligned}$$

and

$$\Gamma_T = \frac{1}{T}\sum_{t=1}^{T} \boldsymbol{x}_t\boldsymbol{x}'_t.$$

Let

$$A = \begin{pmatrix} \Psi_1 \dots & \Psi_q \\ \mathbf{0} & I_{q-1} \end{pmatrix}.$$

Then A has all eigenvalues inside the unit circle and $\boldsymbol{x}_t = A\boldsymbol{x}_{t-1} + \boldsymbol{\varepsilon}_t$. Therefore, there exist a non-singular matrix $\Gamma$ such that

$$\lim_{T\to\infty} \Gamma_T = \Gamma \quad a.s.$$

Furthermore, by Lemma 1 and 2,

$$\lim_{T\to\infty} a_T = 0,$$
$$\lim_{T\to\infty} \boldsymbol{c}_T = \mathbf{0} \quad a.s.$$

Consequently,

$$\lim_{T\to\infty} G_T \sum_{t=1}^{T} \boldsymbol{z}_t \boldsymbol{z}_t' G_T = \begin{pmatrix} 1 & 0 & \mathbf{0}' \\ 0 & 1 & \mathbf{0}' \\ \mathbf{0} & \mathbf{0} & \Gamma \end{pmatrix}$$

Since $\Gamma$ is nonsingular, (18) is proved if

$$\begin{aligned} \log\det(G_T^{-2}) &= \log(\sum_{t=1}^{T} u_t^2) + \log(\sum_{t=1}^{T} v_t^2) + q\log T \\ &\sim (p+2)\log T \quad a.s. \end{aligned} \tag{19}$$

By (4) and (5) of Lemma 1,

$$\lim_{T\to\infty} \frac{1}{\log T} \sum_{t=1}^{T} u_t^2 = 2 \quad a.s.$$

Similar result holds for $\{v_t\}$. Therefore,

$$\log\det(G_T^{-2}) \sim (4+q)\log T = (p+2)\log T \quad a.s.$$

This completes our proof. □

**Remark 1.** Theorem 3 of Wei [17] shows that under similar assumptions as in our analysis,

$$C_T \sim \sigma^2 \log\det(\sum_{t=1}^{T} \boldsymbol{y}_t \boldsymbol{y}_t') \quad a.s. \tag{20}$$

Thus,

$$C_T \sim (p+2)\sigma^2 \log(T) \quad a.s.$$

**Remark 2.** Theorem 1 and Remark 1 have an immediate implication for model selection and can greatly simplify the proof of Theorem 3.5 of Wei [18]. Let $p^*$ be known and $p_0 = \max\{j : \beta_j \neq 0, 1 \leq j \leq p^*\}$ as in (1). Denote $PLS_T(p) = \sum_{t=t_0}^{T} (y_t - \hat{y}_t)^2$ where $\hat{y}_t$ is the forecast of $y_t$ based upon information up to t-1 using the AR($p$) model as in (1) and $PLS_T(\hat{p}_T) = \inf\{PLS_T(j) : 0 \leq j \leq p^*\}$. Wei [18] showed that for both cases of underspecifying and overspecifying AR order ($j$), $\mathrm{P}(PLS_T(j) > PLS_T(p_0) \quad eventually) = 1$. Thus, $\mathrm{P}[\hat{p}_T = p_0 \quad eventually] = 1$. For the case of overspecification, Wei decomposed $\phi_p(z)$ into a sum of a unit root component and a stable component, and worked out the differnece of $C_T$ between the true and the overspecified models. Our results can greatly simplify the proof. Let $C_T^{(j)} = \sum_{t=1}^{T}(y_t - \hat{y}_t^{(j)} - \epsilon_t)^2$ where $\hat{y}_t^{(j)}$ is the forecast of $y_t$ at $t-1$ using the AR($j$) model. For the case of overspecification, $\beta_j = 0, \forall j > p_0$. Applying Theorem 1 and Remark 1, $C_T^{(j)} \to (j+2)\sigma^2 \log(T) > (p_0+2)\sigma^2 \log(T) = C_T^{(p_0)} \quad a.s.$ As for the case of underspecification, $l < p_0$, the desired result, $\mathrm{P}[PLS_T(l) > PLS_T(p_0) \quad eventually] = 1$, is a direct consequence of Theorem 3.2 of Wei [18] since $\beta_{p_0} \neq 0$. Thus, $\mathrm{P}[\hat{p}_T = p_0 \quad eventually] = 1$.

## 3. Implications and applications of the main theorem

We have just proved that for an AR(p) process, $C_T = p\sigma^2 \log(T)$ if it is stationary and $C_T = (p+1)\sigma^2 \log(T)$ if there is an root of 1. Our theorem implies that if the existence of unit root is properly detected and unit root constraint is imposed in forming the forecast, then $C_T = (p-1)\sigma^2 \log(T)$. That is, for model with unit root, estimation is done for the differenced series rather than level of the series. By so doing, we reduce $C_T$ by $2\sigma^2 \log(T)$ which could be substantial for large $T$ and $\sigma^2$. However, it should be noted that $\sum_{t=1}^{T} (y_t - \hat{y}_t)^2$ is not severely affected by existence of unit root since $C_T$, which is of the order of $\log(T)$, is dominated by $\sum_{t=1}^{T} \epsilon_t^2$, which is of the order $T$. This result is natural since it is the long term forecast and not the short term forecast that unit root has strong impact. These findings are further confirmed in our simulation study in Section 5.

In addition, our theorem implies that for AR(p) processes with root equal to or less than 1 in magnitude, as $T \longrightarrow \infty$,

$$\log\det \frac{1}{\log(T)}(\sum_{t=1}^{T} \boldsymbol{y}_t \boldsymbol{y}_t') \longrightarrow c \quad a.s. \tag{21}$$

where $c = (p+1)$ if there is a root of 1 and $c = p$ if all roots are less than one. Equivalently,

$$\hat{d}_T = [\frac{1}{T} \log\det \sum_{t=1}^{T} \boldsymbol{y}_t \boldsymbol{y}_t' - p]^{1/2} \longrightarrow d, \quad a.s. \tag{22}$$

where $d$ is 1 if there is a root of 1 and 0 if there is no unit root. Note that if $p$ is unknown but $r \geq p$ is given, (22) is still true with $r$ replacing $p$ in (22) and in the definition of $\boldsymbol{y}_t$ in (17). In other words, our theorem proves that $\hat{d}_T$ can be used as a test statistic for unit root. This issue will be further investigated in future research.

## 4. Multi-step and adaptive forecast

Our previous analysis focuses on 1-step forecast and there are cases when multiple-step forecast is the main concern. It is conjectured that our results can be extended to multi-step forecast but the issue will be pursued elsewhere. Instead, we shall concentrate our discussion on the relationship between model misspecification and adaptive forecast.

By (1), we have

$$y_{t+h} = \beta_1 y_{t+h-1} + \cdots + \beta_p y_{t+h-p} + \epsilon_{t+h} \tag{23}$$

and

$$\hat{y}_{t+h} = \hat{\beta}_1 \hat{y}_{t+h-1} + \cdots + \hat{\beta}_p \hat{y}_{t+h-p} \tag{24}$$

where $\hat{y}_{t+h-k} = y_t$ for $h \leq k$. So, (24) can be recursively solved in the order of $\hat{y}_{t+1}, \hat{y}_{t+2}, \ldots, \hat{y}_{t+h}$. This is the conventional Box-Jenkins multi-step forecaster.

Another way of generating the multi-step forecast is to solve the model that minimizes the multi-step forecast error and then use it to form multi-step forecast (see Ing [8], Bhansali [2], Weiss [19], and Tiao and Tsay [15]). More specifically, the $h$-step forecast error $e_t(h)$ at time $t$ is

$$e_t(h) = \epsilon_{t+h} + \Psi_1 \epsilon_{t+h-1} + \cdots + \Psi_{h-1} \epsilon_{t+1}$$

where $\Psi_i$ is defined by $[1 - \beta B - \cdots - \beta_p B^p]^{-1} = \Psi_0 + \Psi_1 B + \cdots$. The cost function to be minimized is

$$C(h) = \sum_{t=1}^{T-h} e_t^2(h) \tag{25}$$

Note that for different $h$ different models are used and this explains the name 'adaptive' forecast. Solving (25) involves nonlinear optimization as $\Psi_i$ is a nonlinear function of $(\beta_1, \ldots, \beta_p)$. In practice, approximate linear model is used. That is, the following regression is performed

$$y_t = a_1 y_{t-h} + a_2 y_{t-h-1} + \cdots + a_p y_{t-h-p+1} + b_t$$

and

$$\hat{y}_{t+h} = a_1 y_t + a_2 y_{t-1} + \cdots + a_p y_{t-p+1}$$

The idea behind the adaptive forecast is that if the model is misspecified, that is, p is mistakenly chosen, then this mistake will be amplified radically for the long term forecast. Adaptive forecast could avoid this compounding impact. It is reasonable to expect good performance of Box-Jenkins forecaster for the correctly specified model and good performance of adaptive forecaster for misspecified model.

Ing, Lin and Yu [10] propose a predictor selection criterion to choose the best combination of prediction models (AR lags) and prediction methods (adaptive or plug-in). When there is only one unit root, the proposed method is proved to be asymptotically efficient in the sense that the predictor converges with probability one to the optimal predictor which has minimal loss function.

## 5. Monte Carlo experiments

To assess the theoretical results obtained in previous section and acquire experience about empirical analysis in the sequel, we conduct two Monte Carol experiments. The first is to investigate the finite sample properties of $C_T$ in theorem 1 and the second on forecast comparison between alternative forecasters. For both cases, we generate data from the following four models:

- Model 1: $(1-0.5B)^2(1-B)y_t = \epsilon_t$ or $y_t = 2y_{t-1} - 1.25y_{t-2} + 0.25y_{t-3} + \epsilon_t$. Roots are 0.5, 0.5 and 1.0 respectively.
- Model 2: $(1-0.5B)^2(1-.99B)y_t = \epsilon_t$ or $y_t = 1.99y_{t-1} - 1.24y_{t-2} + 0.2475y_{t-3} + \epsilon_t$
  Roots are 0.5, 0.5 and 0.99 respectively.
- Model 3: $(1-0.5B)^2(1-.95B)y_t = \epsilon_t$ or $y_t = 1.95y_{t-1} - 1.2y_{t-2} + 0.2375y_{t-3} + \epsilon_t$
  Roots are 0.5, 0.5 and 0.95 respectively.
- Model 4: $(1-0.5B)^3 y_t = \epsilon_t$ or $y_t = 1.5y_{t-1} - 0.75y_{t-2} + 0.125y_{t-3} + \epsilon_t$
  All roots are 0.5.

$\sigma^2$ is set to be 1 for all models.

### 5.1. Monte Carlo experiment on $C_T$

The number of replications are 1000 for each experiment. For each, realization, 10 sets of samples are drawn from each model with sample size, T, varying from 100, 200 to 1000. For each sample, starting from $t = t_0(=10)$, the model parameters are estimated and is then used to forecast $t+1$. Then we reestimate the model using sample from 1 to $t+1$ and forecast $t+2$. The process is repeated until when $T-1$ sample is used to estimate the model and then used to forecast $y_T$. The forecast mean square error is then summed from $t_0+1$ to $T$ to obtain $\hat{C}_T$. Finally, we compute the averaged $\hat{C}_T$ obtained from 1000 replications. In other words,

$$\hat{C}_T = \frac{\sum_{i=1}^{1000} \sum_{t=t_0}^{T-1} (\hat{y}_{i,t+1} - y_{i,t+1})^2}{(1000)(T-t_0)} \tag{26}$$

In addition, for each model, we repeat the procedure above with the constraint that one of the root is equal to one. The results are summarized in Table 1. As one can easily see, over 40 millions regressions have to performed to obtain this table and usage of updating formula can significantly reduce the computation burden. In Table 1, the first column is sample size. Results for first model with 0 unit root ($d=0$) and 1 unit root ($d=1$) are put in second and third columns. Results for the other three models are put in columns 4 to 9. Our theory predicts that: (1)

TABLE 1
$C_T$ for simulated data

| | Roots are | | | | | | | |
|---|---|---|---|---|---|---|---|---|
| | 0.5,0.5,1.0 | | 0.5,0.5,0.99 | | 0.5,0.5,0.95 | | 0.5,0.5,0.5 | |
| $T$ | $d=0$ | $d=1$ | $d=0$ | $d=1$ | $d=0$ | $d=1$ | $d=0$ | $d=1$ |
| 100 | 23.47 | 12.33 | 23.47 | 12.33 | 23.80 | 15.36 | 21.07 | 23.22 |
| 200 | 27.55 | 14.71 | 27.55 | 14.71 | 27.71 | 20.19 | 24.28 | 37.60 |
| 300 | 29.90 | 16.06 | 29.90 | 16.06 | 29.83 | 23.90 | 26.09 | 50.86 |
| 400 | 31.49 | 17.00 | 31.49 | 17.00 | 31.21 | 27.17 | 27.32 | 63.57 |
| 500 | 32.75 | 17.73 | 32.75 | 17.73 | 32.26 | 30.27 | 28.29 | 75.96 |
| 600 | 33.76 | 18.30 | 33.76 | 18.30 | 33.09 | 33.12 | 29.04 | 88.12 |
| 700 | 34.62 | 18.79 | 34.62 | 18.79 | 33.80 | 36.01 | 29.69 | 100.41 |
| 800 | 35.38 | 19.22 | 35.38 | 19.22 | 34.40 | 38.89 | 30.26 | 112.72 |
| 900 | 35.99 | 19.60 | 35.99 | 19.60 | 34.94 | 41.65 | 30.76 | 124.80 |
| 1000 | 36.55 | 19.94 | 36.55 | 19.94 | 35.42 | 44.37 | 31.21 | 136.93 |
| $\beta$ | 5.2849 | 2.8583 | 5.2849 | 2.8583 | 5.1947 | 5.2064 | 4.5635 | 13.8536 |
| $R^2$ | 0.9988 | 0.9902 | 0.9988 | 0.9902 | 0.9920 | 0.6315 | 0.9930 | 0.4471 |

Table 2
*MSE for simulated Data*

| | Roots are | | | | | | | |
|---|---|---|---|---|---|---|---|---|
| | 0.5,0.5,1.0 | | 0.5,0.5,0.99 | | 0.5,0.5,0.95 | | 0.5,0.5,0.5 | |
| $T$ | $d=0$ | $d=1$ | $d=0$ | $d=1$ | $d=0$ | $d=1$ | $d=0$ | $d=1$ |
| 100 | 117.81 | 106.92 | 117.81 | 106.92 | 118.33 | 110.06 | 115.59 | 118.17 |
| 200 | 227.10 | 214.61 | 227.10 | 214.61 | 227.36 | 220.24 | 224.11 | 237.66 |
| 300 | 334.88 | 321.53 | 334.88 | 321.53 | 334.97 | 329.57 | 331.53 | 356.59 |
| 400 | 441.30 | 427.33 | 441.30 | 427.33 | 441.27 | 437.73 | 437.66 | 474.10 |
| 500 | 547.43 | 533.00 | 547.43 | 533.00 | 547.19 | 545.69 | 543.55 | 591.33 |
| 600 | 653.42 | 638.61 | 653.42 | 638.61 | 653.01 | 653.67 | 649.36 | 708.93 |
| 700 | 759.51 | 744.24 | 759.51 | 744.24 | 758.92 | 761.63 | 755.18 | 826.46 |
| 800 | 865.20 | 849.45 | 865.20 | 849.45 | 864.35 | 869.13 | 860.55 | 943.18 |
| 900 | 970.92 | 954.95 | 970.92 | 954.95 | 970.01 | 976.96 | 966.16 | 1060.21 |
| 1000 | 1076.76 | 1060.56 | 1076.76 | 1060.56 | 1075.72 | 1084.89 | 1071.91 | 1177.63 |

$\hat{C}_T$ increases linearly with $\log(T-t_0)$ and (2) $\hat{C}_T$ without unit root constraint is 2 times $\hat{C}_T$ with unit root constraint.

We run a simple regression of $\hat{C}_T$ against $\log(T-t_0)$ without intercept for each model and report the regression coefficients and $R^2$ in the last row of Table 1. For column 2 and 3 of the table, the regression coefficients are 5.2849 and 2.8583 respectively while $R^2$ are greater than 0.99 for both cases. In summary, model 1 conforms the theoretical results.

As for model 2, one of the root is 0.99. Since it is the 1-step that is the main concern here, the result is almost the same with model 1. This is consistent with the findings of Lin and Tsay [13] that unit root or not does not matter much for short term forecast.

For model 3, the largest root is 0.95 which is not close to 1 enough. Imposing unit root constraint produces much larger $\hat{C}_T$ and the stable relationship between $\hat{C}_T$ and $\log(T)$ deteriorates greatly as is seen from poor $R^2$. This can be justified by the fact that differencing a stationary process produce a unit root in the MA component which can not be approximated by high order AR. The situation become much worse for model 4 where all roots are equal to 0.5.

For the purpose of comparison, we also report the corresponding conventional MSE ($\sum_{t=1}^{T}(y_t-\hat{y}_t)^2$) for the same 4 models above in Table 2. We observed from the table that contrary to the case for $\hat{C}_T$, the MSE for $d=0$ is about the same as for $d=1$. This confirms our previous analysis that $C_T$, though an important quantity for determining the quality of forecast, is of second order importance as compared to $\sum_{t=t_0+1}^{T}\epsilon_t^2$. For 1-step forecast the distinction between unit root and near unit root does not matter much.

### *5.2. Monte Carlo experiment on short-term and long-term forecast comparison*

This simulation is designed to evaluate the short-term and long-term forecasting performance of alternative forecasters. The number of replications are again 1000. For each replication, 400 observations are generated from the four models above. The first 300 observations are reserved for estimation and then used to produce 1 to 60 steps forecast. Next, the model are re-estimated using the first 301 observations and then used to forecast 1 to 60 steps ahead. The procedure is repeated until when the first 399 observations is used for estimation and the last 1-step ahead forecast is formed. So, we have 100 1-step forecasts, 99 2-step forecasts and 40 60-step forecasts. Then, we compute root mean square error (RMSE) for forecast of each step. Finally, the resulting RMSE is averaged over 1000 replications. More specifically, letting $\epsilon_{i,t}(k)$ be the $k$ period ahead forecast error at time $t$ of the $i$-th replication. Then

$$\text{RMSE}(\ell) = E(\ell) = \sqrt{\frac{\sum_{i=1}^{1000} \sum_{t=300}^{400-\ell} \epsilon_{i,t}^2(\ell)}{(1000)(100-\ell+1)}} \tag{27}$$

The simulation results are put in Tables 3 to 6. In each table, column 1 is steps of forecast, column 2 is the RMSE for model with $p = 3$ and $d = 0$, serving as the benchmark for forecast comparison. Columns 3 to 7 are $E(\ell)$ ratios of model with various $p$ and $d$ to column 2.

From these tables we observe the following. First, for stationary processes, the $E(\ell)$ for the correctly model converges to a constant with the rate of convergence depending upon the value of the root. For root of 0.5, the $E(\ell)$ approach a constant as early as $\ell = 6$ while for root of 0.95 $\ell$ does not stabilize until 30. As for root of .99, it is so close to 1 and $E(\ell)$ is still increasing after $\ell = 60$. For process with unit root $E(\ell)$ increases with $\ell$ for all the whole range of $\ell$. Second, the true model outperforms other misspecified models in forecasting. Third, over-specification of unit results in poor forecast. For the case of model 4 (Table 6) $E(\ell)$ for $d = 1$ is 5% higher than $d = 0$ and jumps to more than 50% for $\ell$ greater than 40. For model 3, one of the root is 0.95 and the forecaster for $d = 1$ is still 45% worse than $d = 0$ though a little better than model 4. As for model 2, one of the root is 0.99 and for up to 20 steps, $d = 1$ fares as well as $d = 0$ and is only 10% worse than the true model at 60-step forecast. Fourth, under-specification of unit root only results in small increase of $E(\ell)$. From column 2 of Table 3, the inefficiency is less than 4% from 1-step to 60-step forecasts. Fifth, under- or over-specification of AR order

Table 3

*Forecasting comparison for simulated data: true $p = 3$, roots are 0.5, 0.5, 1.0*

| Steps | $E(\ell)$ | $E(\ell)$ ratio of MSE to model with $p = 3$, $d = 0$ | | | | |
|---|---|---|---|---|---|---|
| $\ell$: | $p = 3$, $d = 0$ | $p = 3$, $d = 1$ | $p = 2$, $d = 0$ | $p = 2$, $d = 1$ | $p = 4$, $d = 0$ | $p = 4$, $d = 1$ |
| 1 | 3.26 | 99.71 | 103.25 | 102.98 | 100.15 | 99.86 |
| 2 | 7.34 | 99.49 | 102.73 | 102.23 | 100.15 | 99.64 |
| 3 | 11.69 | 99.27 | 102.52 | 101.78 | 100.15 | 99.41 |
| 4 | 15.92 | 99.04 | 102.58 | 101.56 | 100.14 | 99.18 |
| 5 | 19.89 | 98.80 | 102.77 | 101.44 | 100.14 | 98.94 |
| 6 | 23.57 | 98.57 | 103.01 | 101.32 | 100.14 | 98.69 |
| 7 | 26.95 | 98.34 | 103.26 | 101.18 | 100.14 | 98.46 |
| 8 | 30.08 | 98.12 | 103.51 | 101.01 | 100.15 | 98.24 |
| 9 | 33.00 | 97.91 | 103.76 | 100.81 | 100.15 | 98.03 |
| 10 | 35.72 | 97.72 | 104.01 | 100.59 | 100.15 | 97.83 |
| 15 | 47.47 | 97.06 | 105.09 | 99.56 | 100.14 | 97.12 |
| 20 | 57.04 | 96.73 | 105.83 | 98.81 | 100.15 | 96.77 |
| 25 | 65.29 | 96.65 | 106.17 | 98.39 | 100.14 | 96.68 |
| 30 | 72.70 | 96.61 | 106.24 | 98.05 | 100.12 | 96.63 |
| 35 | 79.58 | 96.67 | 106.08 | 97.95 | 100.09 | 96.69 |
| 40 | 85.99 | 96.76 | 105.70 | 97.81 | 100.05 | 96.78 |
| 45 | 92.08 | 96.92 | 105.24 | 97.76 | 100.01 | 96.94 |
| 50 | 98.03 | 97.12 | 104.66 | 97.86 | 99.98 | 97.14 |
| 55 | 103.82 | 97.21 | 104.05 | 97.97 | 99.94 | 97.23 |
| 60 | 109.21 | 97.26 | 103.43 | 97.96 | 99.87 | 97.30 |

TABLE 4

*Forecasting comparison for simulate model: true $p = 3$, roots are 0.5, 0.5, 0.99*

| Steps | $E(\ell)$ | $E(\ell)$ ratio of MSE to model with $p = 3$, $d = 0$ | | | | |
|---|---|---|---|---|---|---|
| $\ell$: | $p = 3$ $d = 0$ | $p = 3$ $d = 1$ | $p = 2$ $d = 0$ | $p = 2$ $d = 1$ | $p = 4$ $d = 0$ | $p = 4$ $d = 1$ |
| 1 | 3.26 | 99.91 | 103.18 | 103.27 | 100.15 | 100.06 |
| 2 | 7.32 | 99.85 | 102.70 | 102.74 | 100.15 | 99.99 |
| 3 | 11.61 | 99.79 | 102.53 | 102.51 | 100.15 | 99.93 |
| 4 | 15.76 | 99.74 | 102.61 | 102.54 | 100.14 | 99.87 |
| 5 | 19.62 | 99.70 | 102.81 | 102.67 | 100.14 | 99.82 |
| 6 | 23.14 | 99.66 | 103.06 | 102.82 | 100.15 | 99.77 |
| 7 | 26.36 | 99.63 | 103.31 | 102.95 | 100.15 | 99.75 |
| 8 | 29.30 | 99.62 | 103.57 | 103.04 | 100.17 | 99.73 |
| 9 | 32.00 | 99.63 | 103.82 | 103.10 | 100.18 | 99.72 |
| 10 | 34.50 | 99.65 | 104.06 | 103.14 | 100.18 | 99.74 |
| 15 | 44.87 | 100.10 | 105.10 | 103.38 | 100.22 | 100.13 |
| 20 | 52.75 | 100.98 | 105.62 | 103.95 | 100.26 | 100.98 |
| 25 | 59.01 | 102.16 | 105.59 | 104.90 | 100.27 | 102.14 |
| 30 | 64.21 | 103.31 | 105.24 | 105.83 | 100.25 | 103.28 |
| 35 | 68.72 | 104.51 | 104.66 | 106.92 | 100.21 | 104.47 |
| 40 | 72.71 | 105.65 | 103.95 | 107.84 | 100.15 | 105.60 |
| 45 | 76.32 | 106.85 | 103.14 | 108.83 | 100.11 | 106.80 |
| 50 | 79.68 | 108.11 | 102.32 | 110.03 | 100.07 | 108.06 |
| 55 | 82.76 | 109.27 | 101.56 | 111.27 | 100.03 | 109.22 |
| 60 | 85.54 | 110.20 | 100.83 | 112.15 | 99.98 | 110.15 |

TABLE 5

*Forecasting comparison for simulated model: true $p = 3$, roots are 0.5, 0.5, 0.95*

| Steps | $E(\ell)$ | $E(\ell)$ ratio of MSE to model with $p = 3$, $d = 0$ | | | | |
|---|---|---|---|---|---|---|
| $\ell$: | $p = 3$ $d = 0$ | $p = 3$ $d = 1$ | $p = 2$ $d = 0$ | $p = 2$ $d = 1$ | $p = 4$ $d = 0$ | $p = 4$ $d = 1$ |
| 1 | 3.26 | 100.87 | 102.92 | 104.61 | 100.16 | 101.00 |
| 2 | 7.19 | 101.59 | 102.50 | 105.05 | 100.17 | 101.70 |
| 3 | 11.20 | 102.35 | 102.39 | 105.87 | 100.17 | 102.43 |
| 4 | 14.93 | 103.15 | 102.50 | 107.02 | 100.18 | 103.22 |
| 5 | 18.24 | 104.01 | 102.69 | 108.34 | 100.19 | 104.06 |
| 6 | 21.11 | 104.92 | 102.91 | 109.71 | 100.20 | 104.94 |
| 7 | 23.59 | 105.87 | 103.11 | 111.09 | 100.22 | 105.86 |
| 8 | 25.72 | 106.85 | 103.29 | 112.41 | 100.24 | 106.81 |
| 9 | 27.57 | 107.85 | 103.46 | 113.70 | 100.26 | 107.79 |
| 10 | 29.17 | 108.88 | 103.60 | 114.95 | 100.28 | 108.79 |
| 15 | 34.68 | 114.26 | 103.91 | 120.95 | 100.35 | 114.08 |
| 20 | 37.59 | 119.71 | 103.34 | 126.69 | 100.37 | 119.46 |
| 25 | 39.16 | 124.73 | 102.28 | 132.01 | 100.32 | 124.42 |
| 30 | 40.04 | 128.94 | 101.30 | 136.47 | 100.26 | 128.60 |
| 35 | 40.61 | 132.40 | 100.65 | 140.17 | 100.19 | 132.01 |
| 40 | 41.02 | 135.16 | 100.27 | 142.89 | 100.14 | 134.75 |
| 45 | 41.31 | 138.01 | 99.99 | 145.70 | 100.11 | 137.58 |
| 50 | 41.50 | 140.97 | 99.79 | 148.84 | 100.08 | 140.50 |
| 55 | 41.61 | 143.76 | 99.68 | 152.10 | 100.06 | 143.26 |
| 60 | 41.68 | 145.57 | 99.66 | 154.16 | 100.04 | 145.04 |

TABLE 6
*Forecasting comparison for simulated data: true $p = 3$, roots are all 0.5*

| Steps | $E(\ell)$ | $E(\ell)$ ratio of MSE to model with $p = 3$, $d = 0$ | | | | |
|---|---|---|---|---|---|---|
| $\ell$: | $p = 3$ $d = 0$ | $p = 3$ $d = 1$ | $p = 2$ $d = 0$ | $p = 2$ $d = 1$ | $p = 4$ $d = 0$ | $p = 4$ $d = 1$ |
| 1 | 3.26 | 105.23 | 100.68 | 108.65 | 100.13 | 104.82 |
| 2 | 5.90 | 110.01 | 100.58 | 114.55 | 100.13 | 109.11 |
| 3 | 7.70 | 115.20 | 100.65 | 121.43 | 100.13 | 113.78 |
| 4 | 8.74 | 120.67 | 100.72 | 128.63 | 100.13 | 118.68 |
| 5 | 9.28 | 126.00 | 100.77 | 135.33 | 100.14 | 123.50 |
| 6 | 9.53 | 130.73 | 100.78 | 140.98 | 100.14 | 127.93 |
| 7 | 9.64 | 134.52 | 100.74 | 145.38 | 100.14 | 131.62 |
| 8 | 9.68 | 137.31 | 100.64 | 148.61 | 100.13 | 134.39 |
| 9 | 9.69 | 139.30 | 100.49 | 150.91 | 100.12 | 136.32 |
| 10 | 9.69 | 140.71 | 100.32 | 152.55 | 100.10 | 137.64 |
| 15 | 9.68 | 144.26 | 99.98 | 156.66 | 100.04 | 140.87 |
| 20 | 9.68 | 145.69 | 99.98 | 158.33 | 100.02 | 142.22 |
| 25 | 9.66 | 147.22 | 99.98 | 160.24 | 100.00 | 143.62 |
| 30 | 9.66 | 148.04 | 99.98 | 161.31 | 100.00 | 144.37 |
| 35 | 9.67 | 148.70 | 99.98 | 162.52 | 100.00 | 144.86 |
| 40 | 9.68 | 148.98 | 99.99 | 162.92 | 100.00 | 145.11 |
| 45 | 9.68 | 150.37 | 99.99 | 164.35 | 100.00 | 146.48 |
| 50 | 9.69 | 153.41 | 99.98 | 167.67 | 99.99 | 149.35 |
| 55 | 9.68 | 155.18 | 99.99 | 169.99 | 99.99 | 150.98 |
| 60 | 9.66 | 155.33 | 99.99 | 170.78 | 99.99 | 150.92 |

only affects the forecast precision marginally. The $E(\ell)$ for all models are within 6% to the true model for all forecasts up to 60-step ahead.

To sum up, the simulation show that slight misspecification of AR order and under specification of unit root are not serious in forecasting but over-specification of unit root could result in poor forecast when the root of characteristic polynomial is far from 1. Yet, improvement of forecasting precision in absolute term could be substantial for large sample when the existence of unit root is appropriately taken into consideration.

## 6. Empirical results

### *6.1. Data*

For empirical analysis, we analyze 6 most frequently used data sets in Taiwan including Gross Domestic Product (GDP), Consumer Price Indices (CPI), Wholesale Price Indices (WPI), Interest Rates( IR), Exchange Rate of New Taiwan Dollar to US Dollar(RX) and money supply(M1B). All series are quarterly data taken from the AREMOS databank. The sample period is 1961:1 to 1995:4 except for M1B which ranges between 1961:3 to 1995:4. So, sample size is 138 for M1B and 140 for the rest series. All series are seasonally unadjusted.

### *6.2. Order selection*

Selecting lag order $p$ and forecasting method simultaneously is analyzed in Ing, Lin and Yu [10]. Here, we follow the conventional wisdom by using AIC and chi-square statistics to determine $p$. When the AIC has a clear minimal, we select the order corresponding to the minimal AIC. When AIC is decreasing without a clear minimum, we use chi-square statistics to select the last significant lag. It turns out that CPI, WPI and RX have order 2, interest rate has order 6, M1B has order 3 and GDP has order 8. The high order indicates the possible existence of seasonal unit root which is not investigated here.

### *6.3. Forecasting procedure*

For each series, the first 100 observations are reserved for estimation and 1- to 20-step forecasts are computed. Then the model are re-estimated using first 101 observations and another 1- to 20-step forecasts are computed. The procedure is repeated until when the first $T-1$ observations are used to estimate the model and the last 1-step forecast is computed. Hence, we have 40 1-step forecasts, 39 2-step forecasts and 20 20-step forecasts except for M1B where there are 38 1-step forecasts and 18 20-step forecasts. For each step, the average root mean square error is computed.

### *6.4. Results*

The results are reported in Tables 7 to 12. From the tables we observe the following. First, $E(\ell)$ increases linearly with $\ell$ for all series except for Interest Rates. This seems to suggest that except IR, all variables have a unit root. Second, regarding the Box-Jenkins forecast, imposing unit root constraint result in poor forecast for all steps ahead for WPI, CPI, GDP and IR. Especially for IR, the RMSE for $d=1$ is 200% higher than that for $d=0$. This seems to be consistent with the finding that its $E(\ell)$ converges to a constant very quickly. However, for RX forecast with $d=1$ fares much better than forecast with $d=0$. The precision gain from imposing unit root is about 5% for 1-step forecast and then up to over 30% for 20-step forecast. This seems to indirectly support the efficient market hypothesis for the foreign exchange market in Taiwan. As for M1B, imposing unit root constraint improves forecast precision from 1- to 7-step forecasts but deteriorates forecast precision from 8-step to 20-step forecasts. The inefficiency is more than 100% for 19 and 20-step forecasts. Third, the performance of adaptive forecaster is mixed. For RX and M1B, adaptive forecast with $d = 0$ and $d = 1$ consistently outperforms conventional Box-Jenkins' forecast by a large margin. The precision gain could go as high as 50%. For CPI adaptive forecast performs poorly for $d = 0$ but very well for $d = 1$. For IR and WPI adaptive forecast with $d = 0$ performs well in short and medium term forecast but fares poorly in long term forecast. But adaptive forecast with $d = 1$ performs okay in the short term but very poorly in the long term. The case GDP is quite interesting. While adaptive forecast with $d = 0$ fares poorly for short and long term forecast, the performance of adaptive forecast with $d = 1$ jumps up and down across steps. This seems to suggest that seasonality plays an important for the differenced GDP which is supported by the corresponding autocorrelation function. This issue will be investigated in future study.

TABLE 7
*Forecasting comparison for GDP*

| | | $E(\ell)$ ratio to model BJ, $d=0$ | | |
|---|---|---|---|---|
| $\ell$ | $E(\ell)$ | BJ, $d=1$ | Adap, $d=0$ | Adap, $d=1$ |
| 1 | 12258.04 | 101.94 | 100.00 | 99.22 |
| 2 | 18024.46 | 104.69 | 101.47 | 179.52 |
| 3 | 21232.71 | 106.87 | 115.10 | 151.59 |
| 4 | 24719.28 | 108.70 | 106.82 | 83.76 |
| 5 | 31938.87 | 110.92 | 102.59 | 99.66 |
| 6 | 37316.81 | 112.05 | 116.19 | 125.65 |
| 7 | 40063.38 | 112.71 | 132.45 | 98.06 |
| 8 | 40505.82 | 109.83 | 147.52 | 45.57 |
| 9 | 46605.77 | 111.94 | 158.02 | 65.51 |
| 10 | 52966.58 | 116.89 | 159.53 | 103.16 |
| 11 | 57551.40 | 121.35 | 157.04 | 91.13 |
| 12 | 59480.32 | 120.18 | 169.46 | 66.06 |
| 13 | 66651.35 | 120.95 | 181.21 | 85.92 |
| 14 | 74760.98 | 123.52 | 184.47 | 109.41 |
| 15 | 79555.22 | 124.66 | 174.53 | 94.18 |
| 16 | 81162.26 | 120.64 | 173.44 | 74.78 |
| 17 | 91194.77 | 117.99 | 185.95 | 61.83 |
| 18 | 99975.45 | 122.09 | 181.84 | 99.95 |
| 19 | 105256.88 | 125.83 | 162.48 | 83.94 |
| 20 | 108809.60 | 122.63 | 162.14 | 55.65 |

TABLE 8
*Forecasting comparison for CPI*

| steps | | $E(\ell)$ ratio to model BJ, $d = 0$ | | |
|---|---|---|---|---|
| $\ell$ | $E(\ell)$ | BJ, $d = 1$ | Adap, $d = 0$ | Adap, $d = 1$ |
| 1 | 1.12 | 99.39 | 100.00 | 101.46 |
| 2 | 1.71 | 98.81 | 83.53 | 77.07 |
| 3 | 1.82 | 99.20 | 100.10 | 69.33 |
| 4 | 1.83 | 99.75 | 123.42 | 79.90 |
| 5 | 2.09 | 98.49 | 128.11 | 84.27 |
| 6 | 2.51 | 99.73 | 124.42 | 81.85 |
| 7 | 2.54 | 101.83 | 148.32 | 89.80 |
| 8 | 2.65 | 102.66 | 165.72 | 95.17 |
| 9 | 2.97 | 102.13 | 163.20 | 95.48 |
| 10 | 3.19 | 101.81 | 171.06 | 95.60 |
| 11 | 3.06 | 106.53 | 209.48 | 98.15 |
| 12 | 3.12 | 108.45 | 237.77 | 96.83 |
| 13 | 3.55 | 106.61 | 236.77 | 86.76 |
| 14 | 3.71 | 107.33 | 257.23 | 89.74 |
| 15 | 3.76 | 111.25 | 289.99 | 98.66 |
| 16 | 3.83 | 113.36 | 326.41 | 99.02 |
| 17 | 4.32 | 110.28 | 336.41 | 91.27 |
| 18 | 4.50 | 111.09 | 372.91 | 89.71 |
| 19 | 4.44 | 113.86 | 431.77 | 86.81 |
| 20 | 4.79 | 111.22 | 455.83 | 80.87 |

TABLE 9
*Forecasting comparison for WPI*

| | | $E(\ell)$ ratio to model BJ, $d = 0$ | | |
|---|---|---|---|---|
| $\ell$ | $E(\ell)$ | BJ, $d = 1$ | Adap, $d = 0$ | Adap, $d = 1$ |
| 1 | 1.16 | 102.47 | 100.00 | 101.19 |
| 2 | 2.13 | 103.43 | 58.78 | 93.97 |
| 3 | 3.05 | 104.54 | 50.13 | 97.35 |
| 4 | 3.88 | 105.55 | 47.53 | 96.59 |
| 5 | 4.56 | 107.43 | 49.39 | 108.12 |
| 6 | 5.07 | 109.64 | 53.84 | 113.24 |
| 7 | 5.41 | 112.38 | 61.27 | 113.87 |
| 8 | 5.58 | 116.54 | 69.88 | 120.83 |
| 9 | 5.89 | 120.00 | 76.71 | 132.83 |
| 10 | 6.36 | 121.66 | 81.36 | 140.82 |
| 11 | 6.93 | 122.71 | 86.87 | 137.69 |
| 12 | 7.53 | 123.95 | 90.96 | 135.62 |
| 13 | 8.02 | 125.92 | 96.49 | 141.24 |
| 14 | 8.48 | 127.93 | 103.98 | 150.77 |
| 15 | 8.82 | 130.43 | 113.97 | 157.73 |
| 16 | 8.87 | 135.06 | 127.60 | 157.82 |
| 17 | 8.97 | 139.02 | 141.47 | 173.71 |
| 18 | 9.15 | 142.11 | 157.87 | 191.82 |
| 19 | 9.52 | 143.65 | 174.20 | 193.95 |
| 20 | 10.19 | 142.43 | 186.53 | 178.22 |

Table 10
*Forecasting comparison for RX*

| $\ell$ | $E(\ell)$ | $E(\ell)$ ratio to model BJ, $d=0$ | | |
|---|---|---|---|---|
| | | BJ, $d=1$ | Adap, $d=0$ | Adap, $d=1$ |
| 1 | .66 | 95.38 | 100.00 | 99.40 |
| 2 | 1.32 | 92.38 | 53.40 | 90.39 |
| 3 | 2.01 | 89.53 | 40.91 | 81.93 |
| 4 | 2.77 | 88.16 | 37.25 | 75.67 |
| 5 | 3.41 | 87.09 | 38.60 | 73.91 |
| 6 | 3.99 | 86.26 | 42.45 | 72.13 |
| 7 | 4.42 | 85.01 | 47.25 | 72.20 |
| 8 | 4.71 | 83.30 | 54.59 | 70.22 |
| 9 | 5.03 | 82.30 | 61.77 | 68.81 |
| 10 | 5.29 | 81.93 | 68.35 | 72.37 |
| 11 | 5.58 | 81.90 | 74.77 | 75.55 |
| 12 | 5.94 | 81.98 | 79.20 | 77.42 |
| 13 | 6.31 | 81.57 | 81.79 | 76.30 |
| 14 | 6.67 | 80.44 | 83.67 | 74.67 |
| 15 | 6.93 | 78.33 | 85.17 | 74.42 |
| 16 | 7.11 | 75.67 | 86.40 | 70.66 |
| 17 | 7.34 | 72.78 | 86.34 | 62.28 |
| 18 | 7.56 | 70.48 | 86.23 | 58.96 |
| 19 | 7.87 | 69.58 | 85.50 | 59.46 |
| 20 | 8.24 | 69.51 | 83.95 | 51.02 |

Table 11
*Forecasting comparison for M1B*

| $\ell$ | $E(\ell)$ | $E(\ell)$ ratio to model BJ, $d=0$ | | |
|---|---|---|---|---|
| | | BJ, $d=1$ | Adap, $d=0$ | Adap, $d=1$ |
| 1 | 96297.39 | 92.81 | 100.00 | 102.13 |
| 2 | 152389.01 | 94.00 | 65.10 | 86.82 |
| 3 | 208305.92 | 94.22 | 52.25 | 78.01 |
| 4 | 266876.68 | 94.73 | 52.56 | 71.28 |
| 5 | 377584.71 | 92.59 | 55.08 | 61.12 |
| 6 | 481271.98 | 94.39 | 61.18 | 56.24 |
| 7 | 532886.35 | 99.13 | 59.46 | 55.99 |
| 8 | 595646.12 | 101.27 | 71.45 | 49.52 |
| 9 | 668073.88 | 108.19 | 84.69 | 50.29 |
| 10 | 774390.51 | 113.50 | 92.29 | 49.31 |
| 11 | 821482.40 | 123.40 | 89.06 | 50.77 |
| 12 | 886619.83 | 129.37 | 90.62 | 46.46 |
| 13 | 1052170.67 | 129.01 | 89.81 | 42.85 |
| 14 | 1158059.45 | 143.04 | 86.93 | 44.22 |
| 15 | 1335812.44 | 145.93 | 70.46 | 44.98 |
| 16 | 1378939.42 | 165.55 | 60.58 | 42.82 |
| 17 | 1649465.49 | 162.79 | 56.17 | 46.78 |
| 18 | 1748042.18 | 183.13 | 61.66 | 45.31 |
| 19 | 1893323.08 | 200.50 | 51.17 | 41.98 |
| 20 | 2079603.50 | 214.50 | 50.45 | 28.34 |

TABLE 12
*Forecasting comparison for IR*

| | | $E(\ell)$ ratio to model BJ, $d = 0$ | | |
|---|---|---|---|---|
| $\ell$ | $E(\ell)$ | BJ, $d = 1$ | Adap, $d = 0$ | Adap, $d = 1$ |
| 1 | 0.72 | 104.58 | 100.00 | 106.60 |
| 2 | 1.16 | 108.68 | 67.94 | 79.35 |
| 3 | 1.24 | 114.84 | 64.70 | 108.65 |
| 4 | 1.38 | 120.28 | 58.85 | 127.00 |
| 5 | 1.63 | 124.48 | 48.70 | 125.50 |
| 6 | 1.82 | 130.14 | 46.05 | 120.60 |
| 7 | 1.83 | 138.25 | 50.44 | 124.52 |
| 8 | 1.82 | 146.64 | 55.74 | 139.96 |
| 9 | 1.83 | 154.97 | 58.67 | 170.56 |
| 10 | 1.83 | 164.17 | 59.64 | 191.38 |
| 11 | 1.76 | 175.17 | 60.49 | 208.55 |
| 12 | 1.70 | 185.17 | 63.08 | 211.01 |
| 13 | 1.67 | 194.23 | 65.27 | 208.26 |
| 14 | 1.61 | 204.93 | 75.07 | 213.25 |
| 15 | 1.52 | 216.39 | 91.16 | 219.55 |
| 16 | 1.36 | 237.94 | 94.70 | 238.53 |
| 17 | 1.27 | 254.73 | 108.89 | 239.50 |
| 18 | 1.10 | 289.28 | 101.44 | 237.00 |
| 19 | 0.84 | 370.74 | 131.31 | 269.90 |
| 20 | 0.59 | 521.61 | 199.32 | 349.14 |

To sum up, the empirical findings are mixed. Imposing unit root constraint might improve forecast precision for some cases but deteriorate forecast precision in others. Also, adaptive forecast differs from Box-Jenkins' forecast by the big margin. Most frequently, it could improve short to medium term forecast but result in poor long term forecast. However, for some cases, it could produce either better or worse forecast for forecast of all steps. Further study is needed to determine the influencing factors.

## 7. Conclusions

We have analyzed the least square forecaster from various aspects. From the theoretical viewpoint, we prove that $C_T$, the most important quantity when evaluating the performance of 1-step forecasters is equal to $(p + d)\sigma^2 \log(T)$ where $d$ is 1 or 0 depending if there is a unit root. This result could be used to analyze the gain in forecasting precision when unit root is detected and is taken into account. Further, this theorem can lead to a simple proof of the strong consistency of PLS in AR model selection and a new test of unit root.

Our simulation analysis confirms the theoretical results. In addition, we also learn that while mis-specification of AR order has marginal impact on forecasting precision over-specification of unit root strongly deteriorate the quality of long term forecast. As for the empirical study using Taiwanese data, the result is mixed. Adaptive forecast and imposing unit root improves forecast precision for some cases but deteriorates forecasting precision for other cases.

## Acknowledgments

Financial support from the National Science Council under grant NSC85-2415-H001-009 is gratefully acknowledged. We would like to thank C. K. Ing and two anonymous referees for helpful comments and suggestions. Without mentioning, the authors are responsible for any remaining error.